\documentclass{amsart}
\usepackage{graphicx}
\usepackage{url}

\theoremstyle{definition}

\theoremstyle{remark}

\numberwithin{equation}{section}

\begin{document}

%\title[short text for running head]{full title}
\title[AlephZero and Mathematical Experience]{AlephZero and Mathematical Experience}

\author{Simon DeDeo}
\address{Department of Social \& Decision Sciences, Carnegie Mellon University, Pittsburgh PA 15123 USA \& the Santa Fe Institute, Santa Fe NM 87501 USA}
\email{sdedeo@andrew.cmu.edu}
\thanks{Contribution to a special issue of the Bulletin of the American Mathematical Society, ``Will machines change mathematics?''. I thank Cris Moore, David Kinney, Jeremy Avigad, and John Bova for helpful discussions. This work was supported in part by the Survival and Flourishing Fund.}

%\subjclass[2020]{Primary }

%\date{\today}

\maketitle

\begin{abstract}
This essay explores the impact of automated proof construction on three key areas of mathematical cognition: on how we judge the role one piece of mathematics plays in another, on how we make mistakes in reasoning about mathematical objects, and on how we understand what our theorems are truly about. It concludes by speculating on a new form of mathematical experience that these methods could make possible: ``glitching'', a game-like search for uncanny consequences of our definitions.
\end{abstract}

\section{Introduction}

\noindent The advent of proof assistants such as \emph{Lean} and \emph{Coq}, combined with progress in Large Language Models and self-play systems such as AlphaZero, raises the question of what happens, to the practice of mathematics, when they are combined. In a recent paper that formed the basis of a 2022 Fields Institute Symposium~\cite{akshay}, Akshay Venkatesh even asks us imagine an ``AlephZero'', trained not on the rules of Go, but on those of mathematical deduction, and that gains, in turn, human or post-human capacities, and is integrated into the mathematical community.

This essay draws on basic ideas in cognitive science to predict three consequences of the contemporary turn to automated methods. First, it predicts a shift in mathematical judgement, as automation eliminates or blurs out experiences of impasse, crucial to anchoring judgements of value. Second, a shift in how we explore the nature of mathematical objects, as automated systems make it harder to state, even temporarily, false things about them. Third, a shift in how we relate to mathematical ideas, as our ability to create truths outpaces our ability to know what they might be about. 

I present these three predictions in turn; my goal is show  how machine technologies might bring them about. To bring the import of these changes into relief, I sometimes talk in terms of loss: to the extent that we rely on automation in certain ways, we will no longer have certain kinds of mathematical experiences associated with impasse, error, and aboutness.\footnote{I use ``aboutness'' here in the philosophical sense of ``intentionality'': the way in which one thing (\emph{e.g.}, a logical formula) can be represent, be about, another (\emph{e.g.}, a shared idea, a personal intuition, or a mental model of a mathematical object).} Some of these experiences will be fewer in number, or lower in resolution---vaguer, briefer, less precise. Parts of my discussion might, as a consequence, remind readers of a dystopian future by the writer Ted Chiang~\cite{chiang2000catching}, where humans are outcompeted by oracular ``metahumans'', and cease seeking knowledge altogether.

It can be useful to imagine such futures, because they can serve as intuition pumps~\cite{dennett2013intuition}. However, the nature of human curiosity suggests that now is not the time to expect them. Impasse, error, and aboutness will not disappear, and mathematicians at the cutting edge of automated methods (\emph{e.g.}, Patrick Massot~\cite{massot}) provide vivid accounts of all three. Those same accounts, however, also emphasise the ways in which their experiences are fundamentally different from what has come before. In the final section, ``Glitching, Clipping, and Logical Exploits'', I speculate on the future evolution of this process.

\section{Value and Impasse} %Testimony and Experience}
\label{value}

Venkatesh's account of mathematical value emphasizes the importance of a conjecture being \emph{central}: a conjecture acquires value when it is ``linked with many other questions of (prior) importance''~\cite{akshay}; Jeremy Avigad~\cite{avigad2022varieties} uses similar language. Mathematicians seem to follow a heuristic familiar to both the sciences and day-to-day reasoning: just as we often value, even over-value, explanations that link together an apparent diversity of prior observations~\cite{zach}, mathematicians value conjectures that link together prior questions.

On the surface, such a criterion is both intuitive and clear. Consider the Langlands Program, a celebrated example of value-through-linkage in modern mathematics. Even I, a non-mathematician, have heard of Langlands---for example in Michael Harris' admirable \emph{Mathematics without Apologies}~\cite{harris2017mathematics}. What I hear sounds like something that ought to be valued very highly indeed.

If pressed, however, on what the Langlands Program \emph{actually} is, I would say that it is an attempt to prove a set of theorems that link together high-value but hard-to-prove facts in number theory (on the one hand) and geometry (on the other). There are definitely curves involved, and integers, but what characterizes the theorems that make it into the Program and why the particular correspondences they govern, rather than others, are the object of such deep fascination, remains mysterious to me. 

My knowledge of the value of the Langlands Program is partial. This is not (just) because I don't know the theorems implicated in its correspondences, but also because my acquaintance with how this or that correspondence plays out when trying to prove things is at (at best) second-hand. My beliefs about the Program's value-relevant properties (``centrality'') come to me through the testimony of people who have tried to prove things when a Langlands correspondence is relevant, and not through the mathematical experience of trying to do so myself. %I have not spent days, say, staring at a blackboard, stumped, until one of my colleagues realizes that this step, here, can be transformed into a statement that feels analogous to some fragment of of a Langlands correspondence, and I have not followed out her reasoning to see how this might be so. Without this kind of experience, I am forced to rely on the testimony of others who have had experiences like this.

Even if my beliefs about that value are correct, in other words, there is something suspect about my holding them without, at least silently, adding, ``according to those who know''. Intuitively, the situation is similar to the phenomenon of testimony in aesthetic matters~\cite{robson2012aesthetic}: for me to wax enthusiastic about the ``deep importance'' of Langlands, and the ``true centrality'' of its aims, is akin to someone praising the prose style of a book he has never read or the transporting delights of a cathedral he has never visited---with the added twist that even an untutored tourist can experience something of the aesthetic by simply stepping inside the cathedral, while the mathematical case, more demanding, requires some kind of working-through. Adapting Hill's~\cite{hills2022aesthetic} account of aesthetic testimony to the judgement of mathematical value, we might say that Harris' book can transmit the correct beliefs about the value of the Langlands Program, but not the necessary understanding of those values. This can only come from the experience of working with Langlands itself. % \footnote{Say I have spent weeks at a blackboard, stumped on how to prove a theorem I care about; one of my colleagues comes in, and shows me how ``this step, here''---a crux of my argument, that captures the essence of the larger project I have in mind---can be transformed into a statement that feels unexpectedly analogous to some fragment of a Langlands correspondence. I see how it might be so. Experiences like these, of impasse and difficulty, followed by consilience, the recognition of correspondences, and new forms of mystery, seem to be implicated in the mathematical community's evaluation of the value of the program.} 

We need not restrict ourselves to something as exalted as Langlands. Similar concerns apply to any proof or conjecture, and the judgements of the depth to which it connects to the rest of mathematics. The  judgement that Theorem A's involvement in Theorem B is important, depends, in the final analysis, on someone trying to prove Theorem B, and experiencing, directly, both the difficulty of the impasses that A resolves, and the ways in which A resolves them.\footnote{In the usual process, these experiences are shared socially, where they become testimony---not just to outsiders, but also to other mathematicians without the time or training to experience the process directly. There is nothing intrinsically wrong with this: just as in the case of aesthetic judgement, it is not always wrong to rely on testimony, for example, in the awarding of grants and prizes. C. Thi Nguyen~\cite{thi} argues for an asymmetry in the aesthetic case: testimony can establish that something \emph{ought} to be found beautiful upon acquaintance, but not that it is, indeed, so.}

It seems difficult to eliminate the importance, for mathematical judgement, of working something through. In the sciences, by contrast, I can discover the key relationships between propositions simply by varying their likelihoods of being true and seeing how one affects another~\cite{zach}: no matter how complex the underlying theory for why, say, $X$ affects $Y$, I can determine the importance of a binary variable $Z$ by turning it on and off (or, in probabilistic theories, by making it more or less likely to be on). No similar analysis works in the mathematical case, because I don't know how mathematics (or even logic) works in a world where (say) there were only a finite number of primes. 

One can imagine varying a mathematician's confidence in the validity of a step in a proof---this is what my colleague Scott Viteri and I do in~\cite{viteri}---but what this leads to is an (approximate) cognitive account of mathematical experience, not a paraconsistent theory of mathematics itself. If we see mathematics as a formal, ``timeless'' deductive system, this kind of counterfactual reasoning seems hard to get off the ground: we can say at best that, in mathematics, everything depends on everything and every (proven) theorem depends on every other simply because, if any of them were false, mathematics itself would be inconsistent. Mathematicians have a different way of making sense of importance: rather than asking the (somewhat nonsensical) question ``how much would it matter if such-and-such a theorem, known to be true, was actually false'', they can use their experiences of impasse and resolution as guides for which theorems to value.

Automated methods seem to pose a challenge to how these judgements get made, however, because testimony about value no longer necessarily refers back to human experiences of impasse: it is the humans that get stuck, not the machines, who simply grind through a mechanical process for some period of time. The very goal of these tools is to take some part of the reasoning process out of human experience, and into a more reliable, mechanical realm.\footnote{The transposition is---we hope---truth preserving, \emph{i.e.}, we are just as justified, if not more justified, than before in believing in the truth of the result. The argument here concerns testimony about value, not (as in~\cite{Williams1972}) judgements about truth.} 

These methods might tell us that an efficient proof of Theorem B draws, in its abstract syntax tree,\footnote{The abstract syntax tree (AST) is a representation of how a proof, written in a high-level language like \emph{Coq} or \emph{Lean}, is broken down into its constituent deductive steps; as described in Ref.~\cite{viteri}, an AST is a natural way to summarize the logical dependencies between different statements.} on a fragment of Theorem A. By their very nature, however, these trees are extraordinarily complicated, and simply knowing which theorems are cited is not decisive: Theorem B will depend upon many other theorems, including a host of more or less trivial things, and it may also depend on deep things but only in what a human would consider a trivial fashion. (This is the best-case scenario: a truly efficient proof of Theorem B might, in fact, make the role of Theorem A's lemma harder to notice, distributing it in fragments across the entire tree.)

In response to this challenge, we might conduct a post-hoc analysis of the syntax tree of the proof of Theorem B and show how some graph-theoretic property shows that subtrees, associated with Theorem A, are especially ``central''; my colleagues in network science, for example, might use ``betweenness centrality''~\cite{freeman1977set}. What we gain, however, is not the relevant experience of value, but only knowledge about a hypothesised proxy for value, an operationalization. Such a proxy may correlate with value judgements we already believe from experience (``the role of Theorem A is graph-theoretically similar, at $p<0.01$, to the celebrated role of...''). Proxies, however, are not the thing itself: at best, they are another form of testimony, one that lacks even the backing of experience.

Deployed widely enough, the reliance on such proxies---even if they correlated perfectly with ideal judgement---would lead to a strange scenario: a kind of zombie mathematics, where mathematicians celebrate a theorem not for how it untangles and reorders their reasoning, or the reasoning of their colleagues, but because it has a high centrality score.

Such a dystopian fantasy is unlikely, of course, to actually occur. We can expect automated methods to hide, behind mechanical processes, certain experiences of impasse and resolution. Mathematicians, however, will continue to search for new understandings, ideas, and experiences, and we can expect the dynamics of this process to have experiences of impasse and resolution as well. We should, however, expect those experiences to transform; we return to how they might in Section~\ref{clipping}.

\section{Models and Errors}
\label{error}

Mathematicians make errors. What those errors are doing, however, is particularly hard to understand in the formalist picture, whether that be the formalism that motivated Hilbert's project, or later versions such as Haskell Curry's~\cite{curry1951outlines}. If mathematics is  a matter of logical deduction from axioms, then errors are nothing more than ill-formed thoughts---on a par with claims about ``square triangles'' or ``the third even prime'', forms of nonsense that are not truly thoughts about anything at all~\cite{putnam}. Such errors may be of interest to psychologists, perhaps, as physical causes of belief, but have no relevance for the subject of mathematics itself beyond the bare fact of their ungrammaticality.

Putting ``AlephZero'' systems to the side for a moment, automated proof assistants such as \emph{Coq} and \emph{Lean} inherit, by necessity, formalism's dim view of error. Assuming that the underlying code implements the formal system correctly, a proof assistant can not allow a mathematician to introduce a claim that is not justified on the basis of the axioms and rules of the system. Checking that a step in a proof is correct is similar to checking that a Python program is well formed: a purely mechanical process equivalent to tasks like balancing parentheses. This is true not only when an assistant is used to checking a completed proof (\emph{e.g.}, one submitted to a journal), but also in the course of proof construction. 

My own experience with \emph{Lean}'s tutorial problems---simple proofs in Peano arithmetic designed to teach the basics of the system---was instructive; in contrast to pen-and-paper, which will accept all sorts of nonsense, Lean's interface will only accept, step by step, the justifiable. If one works solely with the interface, one's nonsense is detected immediately on entry, an experience that contrasts with more familiar experiences of mathematical derivation which often involve back-tracking to root out falsehoods. For a proof assistant an unjustifiable step is only an endpoint, a dead-end; there can be no felicitous contradictions that persist, generatively, for a time, before being fixed.

The strict exclusion of error is a defining feature of the automated proof assistant, and even those who subscribe to the formalist picture might worry about how it might impose itself onto the human task of exploration. Because of this, it is particularly interesting to see how mathematicians make sense of what happens when they make errors. 

Some accounts are purely psychological, in a trivial sense; the eight errors listed by Scott Aaronson in~\cite{scott}, for example, point to standard human failures such as hubris and the Dunning-Kruger effect. In these cases, error is no more relevant to mathematical experience than sleep is: we have limits on our ability to be reasonable, and these limits can stop us doing mathematics. In accounts like these, the error-exclusion made possible by proof assistants is a relatively uncomplicated good, to be welcomed in the same way we might welcome lane-keeping alerts for long-haul truck drivers.

Other accounts of error are more cognitive, however, and complicate the story because they see error as something that emerges organically from the process of mathematics itself. Terrence Tao~\cite{tao}, for example, describes the difference between ``local'' and ``global'' errors in a proof, and how the former have mathematical properties that make them easier to fix; statistical study of proofs~\cite{viteri} suggests that their logical structure is indeed modular in ways that help make sense of Tao's local-global distinction. In discussions with mathematicians, one learns that they have a variety of heuristics they use to identify errors in their reasoning, including particular signatures and traces that error leaves on downstream steps. One sign of an error, as I learned during the 2022 Symposium, is that subsequent results become ``too easy'' to obtain, in a kind of limited principle of explosion that temporarily levels the hierarchy of value. On the grandest scale, some mathematicians are celebrated for fertile errors; the mathematician Solomon Lefschetz, for example, published a proof of what is now known as the ``Hard Lefschetz Theorem'' in 1924---a proof that was incorrect, but nevertheless generative enough to found a field.\footnote{In the words of one biographer, ``Lefschetz never stated a false theorem nor gave a correct proof''~\cite{nas}.}

Error-making can sometimes even take on a deliberate form. Terrence Tao begins one argument, on the website \emph{Stackoverflow}, with the counterfactual ``Note that if $\pi$ were rational...'',\footnote{\url{https://mathoverflow.net/questions/282259/is-the-series-sum-n-sin-nn-n-convergent/282290#282290}} using this unpromising beginning to structure the subsequent argument. A deeper level of engagement with error is provided by one correspondent (a reviewer for this article) who described, \emph{ad lib}, the following process:
\begin{quote}
One sometimes assumes the existence of a mathematical object that might be wrong, while trying to prove something hard, but you are confident that there will be a modification that works. However that will probably be a tedious lengthy detail, so better to make progress on the hard question by assuming some version of that object exists, doing the more interesting part of the argument with that object and then going back and getting the original object correct, especially in the light of how it is to be used.
\end{quote}
In this case, the mathematician defers. He stipulates the existence of an object he suspects to have unjustifiable (or even logically inconsistent) properties. The construction of this object allows him to continue the argument, albeit in a bracketed (\emph{i.e.}, logically unjustified, potentially contradictory) form. From this he gets the shape of this argument, including a sense of the purpose his object is required to play. This, in turn, allows him to go backwards in the derivation, to repair the errors in the object's definition---errors without which the original progress could not have been made. Important information is gained by being wrong: without making (or at least risking) the error, it is impossible to get the shape of the future argument, which contains the insight necessary to make the correction.

Mathematics is difficult, and so it is hardly unexpected that mathematicians have developed accounts of how errors can be a source of insight. In the process of finding a proof, mathematicians make mistakes all the time, and if mathematics only happened when one was putting one justified formula in front of another, then the majority of mathematicians would be spending the majority of their time speaking nonsense, no mathematics at all, and, as Lear says to Cordelia, ``nothing can come of nothing''.

If, by contrast, we believe we can experience errors of logic as having mathematical meaning---if we are willing to say, for example, that one can assert falsehoods while making progress---we must also grant that some essential aspect of the activity goes beyond the formalist picture enforced by a proof assistant.

In as much as automated methods help us be wrong less often, in other words, they foreclose an aspect of the mathematical experience. Some might wonder, however, if the generative errors that mathematicians rely on today might be adapted to an era of perfect reason. Tao's argument beginning with ``if $\pi$ were rational'', for example, might be easily rephrased---if less intuitively---to begin ``because $\pi$ is not rational''. Is anything truly lost along with the experience of error?

From the cognitive point of view: very possibly, yes. This is because mathematicians, just as much as any other human, are expected to rely in part on the construction of mental models~\cite{craik1967nature, johnson2001mental, avigad2021reliability}; in particular, reduced and partial mental models of the mathematical objects in play. 

In day-to-day life, mental models are tuned by interaction with the world, through a process of learning and feedback; to develop our mental models, we use them to generate predictions about the world that we then compare to reality~\cite{clark_2013}. If our model fails to predict correctly, we update it. This update is far from instantaneous: sustained model failures help direct our attention---we attend more closely to aspects of our experience that our models failed to predict~\cite{itti}, and prediction failures seem to be not only a core feature of low-level cognition~\cite{pred} but useful guides to high-level decision-making processes such as scientific exploration~\cite{murdock2017exploration,foster2021surprise}. The mathematical parallel is, most naturally, the use of mental models to be wrong about mathematics; without the possibility of being wrong, the mental model can not change and the process of attention is diffused. This leads to a second dystopian fantasy, one where automated methods lead to a world in which our mental models become increasingly vague and low-resolution.

Such a loss would be more than simply epistemic. It is not just that mathematicians would have more impoverished representations of what they are doing, but also that they would be deprived of particular experiences. Mathematicians may have a horror of error, but wandering into error is an unavoidable consequence of taking deliberate risks with our mental models---an experience inseparable from the act of exploring the world with curiosity.

Just as in the previous section, I suggest that this is unlikely to happen: mathematicians are simply too curious about their objects not to want to explore and reason about them in ways that allow them to be wrong.  The most obvious way to continue doing that is in trying---and, naturally, sometimes failing---to predict what an automated theorem prover will do given an input. 

This is only a partial compensation, however, because ``cyborg'' errors last only as long as it takes to phrase and type the first step of the erroneous intuition. Automated systems such as \emph{Lean} can defer a subgoal in a proof with a keyword such as {\tt sorry}, but they are not (yet) able to humor their human counterparts by suggesting that a subgoal has been achieved when the assumption is false. Phrased positively, this section suggests a role for machines that give us new powers by allowing us to be wrong in new ways.

\section{Definitions and Aboutness}
\label{definition}

In his talk at the 2022 Fields Institute Symposium about the Liquid Tensor Experiment~\cite{scholze2022liquid} and its challenge to formalize, in \emph{Lean}, a key theorem of Peter Scholze's work in Condensed Mathematics, Johan Commelin draws attention to an illuminating risk of formalization: the risk that, through malice or accident, one's definitions may make a proof trivial. 

A casual observer of the automated proof can look at the \emph{Lean} statement of the main theorem, and see how it parallels, more or less, Peter Scholze's original human text. There are references to profinite sets, $p$-Banach spaces, and so forth, in the expected places. In the end, however, when we glance back and forth between Scholze's LaTeX statement and \emph{Lean}'s monospaced font,
\begin{quote}
...as Magritte told us ``Ceci n'est pas une pipe''. We could have done something very evil\footnote{Commelin's reference to ``evil'' may seem extreme, but it is not entirely out of bounds; consider, for example, the recent (Spring 2023) scandal of the ``Space Zoom'' feature on Samsung phones: unknown to the public, the internal code was able to recognize when a user was taking a picture of the moon, and filled in details of the image that would have been impossible for the detector itself to have captured given atmospheric conditions.} or we could have done something very stupid---what if we had just defined Ext groups to be zero to begin with, because that's what the main statement is about. We need to prove that for all $i$ some Ext group is zero, well, if we just define it to be zero then we're done ... even if we're not evil we could have done something stupid that would have completely trivialized the proof. [Johan Commelin, ``Abstract Formalities'', 2022 Fields Institute Symposium, minute 59]
\end{quote}
The larger concern, as Commelin articulates it, is that automated methods often place us in a position where we have thousands of definitions, and even in the presence of abstraction boundaries, documentation, and test cases, one is, in the final analysis, thrown back on abductive reasoning to figure out what, exactly, the proof is about. 

We might say, in turn, that we have two forms of aboutness in play. On the one hand, we have the way in which Scholze's proof is an attempt to reason about ideas\footnote{I follow Scholze's use of the word ``idea''; see \url{https://bit.ly/scholze_harris}.} that Scholze had in mind, ideas created through a (collaborative, iterative) process of forming intuitions, making definitions that answered to them, and attempting to prove things. On the other hand, we have the way in which a set of code fragments, labelled ``definitions'' and input into \emph{Lean}'s axiomatic system, are mutually constrained by the laws of type theory to output some sets of code fragments downstream, and not others.

Commelin's remarks speak to the task of mediating between these forms: poetically speaking,  asking if the ``shape'' of the code fragments that the \emph{Lean} proof enables (the representation of the pipe) matches the ``shape'' of Scholze's ideas (the pipe itself). This itself is part of the Liquid Tensor Experiment, which addressed the problem by presenting not just the \emph{Lean} formulations of key steps in Scholze's proof, but also a ``test suite'', that enables users to play with the objects that appear in those formulations, to see if they behave in the ways they ought.

The need to square competing forms of aboutness is reminiscent of the dialogues in Lakatos' \emph{Proofs and Refutations}~\cite{lakatos1963proofs}, which provide an extended example of how students, challenged to prove something, go back and forth exchanging claims and counterexamples, gradually realizing the limits of their intuitions---that what they think they are proving things about is not quite what they think it is---and correcting the terms of the proof as they go. Catarina Dutilh Novaes~\cite{catarina} formalizes this idea as a game between ``prover'' and ``skeptic''; the end point---or the asymptote---of this dialectic is the justified belief that the proof (in the end) truly is about the ideas in question (in the end).

It is unclear, however, what it would mean to achieve intersubjective agreement with a machine. It is, certainly, possible for two human mathematicians to come to agreement on what a \emph{Lean} proof is about, but this is a distinct  task---more analogous to two scientists attempting, through experiments, to determine the most efficient account of the causal structure of a black box, and perhaps the ``stochastic mathematical systems'' proposed by Wolpert and Kinney~\cite{wolpert2023stochastic}. Commelin's reference to the need for abductive reasoning about a \emph{Lean} proof, and the test suite than enables it, may point in this direction, but it is important to distinguish this from the students in Lakatos's dialogues. The Lakatos students are engaged in dialectical logic, not empirical modeling: they argue \emph{with}, rather than \emph{about}, each other. It is this ``arguing with'' experience, and the errors and vagueness it involves, that automated methods seem to exclude.

\section{$\aleph(0)$ in the Loop: Glitching, Clipping, and Logical Exploits}
\label{clipping}

Much of what we know about the impact of automated methods comes from proof verification and assistant systems. These somewhat constraining tools, however, are already being combined with systems like OpenAI's ``Codex'' code-completion technology, which can propose next steps, respond dynamically to criticism and coaching, and which maintains hidden states potentially analogous to mental models of the task at hand. Mathematicians who become familiar with the syntax of these tools may soon be engaged in quite extensive forms of co-construction, where the machine not only fills in small gaps, but proposes new definitions, attempts to prove lemmas of its own devising, and even enlists the human partner in new goals that may be only partially transparent.

One analogy to this experience, sometimes made by those on the cutting edge of this process, is to the iterative feedback found in a video game. If the dialectic of traditional mathematics is somewhat like the human-on-human play of Dungeons and Dragons---with the prover serving as ``Dungeon Master'', and the skeptics attempting to play within, or subvert, the boundaries of the prover's vision---the new era heralded by projects such as the Liquid Tensor Experiment might be thought of a massively multiplayer online role-playing game (MMPORG), with an engine behind the screen that responds, in sometimes counterintuitive ways, to the community's inputs.

Once moves in a proof system are allowed to produce side-effects (\emph{i.e.}, once we move beyond the functional programming paradigm of a system such as \emph{Lean} to more general Codex-like systems), the analogy to video games is tighter than it might at first appear. As long as they do not crash, violate memory boundaries, or fall into an infinite loop, proof assistants and video games are just ordinary computer programs, a complex system of inter-related states and state-to-state transition rules. Even within the functional programming paradigm, game-like experiences seem to be common for those who engage at sufficient depth.

Video games are inherently pleasurable creations, and the 21st Century may well be defined by this novel form~\cite{walz2015gameful}. Unlike other forms of art, they give agency a central role~\cite{nguyen2020games}, which parallels key aspects of mathematics: the experience of mathematics is an agentic one, one of dialectic, choice-making, puzzle-solving, backtracking, and error---far removed from the passive contemplation of timeless structures urged by early versions of formalism.

Even if automated systems take on a video game quality, however, we should not expect mathematicians to be content with ``playing'' an automated proof system the way an ordinary person plays, say, Super Mario Bros. In ``Mathematics and the Formal Turn'', for example, Jeremy Avigad notes how some \emph{Lean} users talk in terms of ``golfing a proof'', a practice analogous to the ``code bumming'' described in early histories of the computer revolution~\cite{levy1984hackers}---and reminiscent of the practice today of the ``speed run'', where someone attempts (say) to finish all the levels of Super Mario Bros in an absolutely minimal amount of time.

One phenomenon not yet explored (to my knowledge) is the attempt to discover, and make use of, glitches. Glitches are a particularly uncanny source of interest, which emerge out of how ordinary video games attempt to simulate a physical reality for the user to explore. Misjudgements in the designers' vision sometimes means that particular objects or locations, despite following the prescriptions of the ``physics engine'' perfectly, have exceptional properties that violate our mental models of the intended reality and produce novel game logics. In one game, for example, the chain on a swing set in a playground---a minor background object intended mostly as scenery---serves as a fount of energy that can launch a player into the sky; in another, a player can walk through walls (``clipping'') if he holds a spoon in front of him as he moves.

Glitches arise at the intersection of the technological and the human~\cite{menkman2011glitch}. They are semantic errors (\emph{i.e.}, they reveal that the physics engine in question is not actually mirroring the world the player expects), but not syntactic ones (\emph{i.e.}, the code has, indeed, compiled correctly, and there is no crash, buffer overflow, or subversion). Following the discussion of Section~\ref{definition}, it seems likely that such glitches should appear in automated systems even when---indeed, precisely because---the type-checker is working as expected. They occupy a space between the intuitively true and the logically false---not true antinomies or logical errors, but uncanny experiences of what ought not to be.

Video game glitches are difficult to find by examination of the code. They tend to be discovered by communities engaged in obsessive forms of play well beyond ordinary use~\cite{bainbridge2007creative}, and similar levels of obsessiveness may be needed for the mathematical case. Glitches might also be found by automated means: for example, fuzzing~\cite{fuzz}, a technique for finding security vulnerabilities in code by trying trillions of randomly-chosen, but synatactically valid, inputs. Fuzzing could be applied to a system like \emph{Lean} and checked against basic human intuitions to seek the unexpected in our axioms. 

Imagine, for example, that the Liquid Tensor Experiment continues, spreading to a wider community of mathematician-hackers who wish to go beyond simply replicating Scholze's original proofs. Some might take a more or less traditional route, building mental models of Scholze's objects, asking what else might be true of them, developing new conjectures, and trying to prove them. Others, however, might rely on the affordances of \emph{Lean} to create proofs whose objects and goals do not answer to any natural human need---whose aboutness has detached from the constraints of an idea. This could be done, for example, by search processes that attempted to maximize interesting structural features of a proof's dependency network, deviated in unusual ways from the structures we have seen before, or even by brute force fuzzing.

One possible outcome for this second group is the unproductive: there are an infinite number of mathematical truths, and most of them do not add up to much. It may well be that human intuition creates objects whose most interesting consequences are those turned up by the kind of exploration humans are naturally drawn to. But there are reasons to believe the outcome might be more interesting than a sampling of complicated but trivial facts. The affordances of a proof assistant, combined with the ingenuity of the hacker, could lead to an exploration process that turns up objects with interesting properties unlike anything we might have imagined for ourselves. 

Such objects would be unusual in mathematical culture: created not for the purposes of answering to human intuition, they would possess, nevertheless, a derivative intentionality, a complexity that emerged from the systems we ourselves created. Outpacing our mind-grown ideas, these objects would have a cyborg-like, even somewhat alien character, and it is not hard to imagine us reasoning about them both abductively, building test suites to see how they behaved, and deductively, seeing what we could prove about them, in turn.

The day when we discover unexpected, uncanny glitches in the definitions we have to hand, or the ones we will co-create with our machines, may be nearer than we expect. Bertrand Russell once wrote, partly in jest, that ``mathematics may be defined as the subject in which we never know what we are talking about, nor whether what we are saying is true''~\cite{russell}. With the advent of AlephZero it seems likely that, while we will continue to not-know these things, we will come to not-know them in new and unexpected ways.

\bibliography{main}

\providecommand{\bysame}{\leavevmode\hbox to3em{\hrulefill}\thinspace}
\providecommand{\MR}{\relax\ifhmode\unskip\space\fi MR }
% \MRhref is called by the amsart/book/proc definition of \MR.
\providecommand{\MRhref}[2]{%
  \href{http://www.ams.org/mathscinet-getitem?mr=#1}{#2}
}
\providecommand{\href}[2]{#2}
\begin{thebibliography}{10}

\bibitem{scott}
Scott Aaronson, \emph{Eight signs a claimed {P}$\neq${NP} proof is wrong},
  (2010), \url{https://scottaaronson.blog/?p=458}.

\bibitem{avigad2021reliability}
Jeremy Avigad, \emph{Reliability of mathematical inference}, Synthese
  \textbf{198} (2021), no.~8, 7377--7399.

\bibitem{avigad2022varieties}
\bysame, \emph{Varieties of mathematical understanding}, Bulletin of the
  American Mathematical Society \textbf{59} (2022), no.~1, 99--117.

\bibitem{bainbridge2007creative}
Wilma~Alice Bainbridge and William~Sims Bainbridge, \emph{Creative uses of
  software errors: Glitches and cheats}, Social Science Computer Review
  \textbf{25} (2007), no.~1, 61--77.

\bibitem{chiang2000catching}
Ted Chiang, \emph{Catching crumbs from the table}, Nature \textbf{405} (2000),
  no.~6786, 517--517.

\bibitem{clark_2013}
Andy Clark, \emph{Whatever next? predictive brains, situated agents, and the
  future of cognitive science}, Behavioral and Brain Sciences \textbf{36}
  (2013), no.~3, 181–204.

\bibitem{craik1967nature}
Kenneth James~Williams Craik, \emph{The nature of explanation}, Cambridge
  University Press, 1967.

\bibitem{curry1951outlines}
Haskell~Brooks Curry, \emph{Outlines of a formalist philosophy of mathematics},
  Elsevier, Amsterdam, NL, 1951.

\bibitem{pred}
Hanneke Den~Ouden, Peter Kok, and Floris De~Lange, \emph{How prediction errors
  shape perception, attention, and motivation}, Frontiers in Psychology
  \textbf{3} (2012).

\bibitem{dennett2013intuition}
Daniel~C Dennett, \emph{Intuition pumps and other tools for thinking}, WW
  Norton \& Company, 2013.

\bibitem{foster2021surprise}
Jacob~G Foster, Feng Shi, and James Evans, \emph{Surprise! measuring novelty as
  expectation violation}, SocArXiv (2021).

\bibitem{freeman1977set}
Linton~C Freeman, \emph{A set of measures of centrality based on betweenness},
  Sociometry (1977), 35--41.

\bibitem{fuzz}
Patrice Godefroid, \emph{Fuzzing: Hack, art, and science}, Commun. ACM
  \textbf{63} (2020), no.~2, 70–76.

\bibitem{nas}
Phillip Griffiths, \emph{{Solomon Lefschetz}}, Biographical Memoirs: Volume 61,
  The National Academies Press, 1992.

\bibitem{harris2017mathematics}
Michael Harris, \emph{Mathematics without apologies}, Princeton University
  Press, 2017.

\bibitem{hills2022aesthetic}
Alison Hills, \emph{Aesthetic testimony, understanding and virtue}, No{\^u}s
  \textbf{56} (2022), no.~1, 21--39.

\bibitem{itti}
Laurent Itti and Pierre Baldi, \emph{Bayesian surprise attracts human
  attention}, Advances in Neural Information Processing Systems \textbf{18}
  (2005).

\bibitem{johnson2001mental}
Philip~N Johnson-Laird, \emph{Mental models and deduction}, Trends in Cognitive
  Sciences \textbf{5} (2001), no.~10, 434--442.

\bibitem{lakatos1963proofs}
Imre Lakatos, \emph{Proofs and refutations}, Nelson London, 1963.

\bibitem{levy1984hackers}
Steven Levy, \emph{Hackers: Heroes of the computer revolution}, Anchor
  Press/Doubleday Garden City, NY, 1984.

\bibitem{massot}
Patrick Massot, \emph{Why formalize mathematics?},
  \url{https://www.imo.universite-paris-saclay.fr/~patrick.massot/files/exposition/why_formalize.pdf}.

\bibitem{menkman2011glitch}
Rosa Menkman, \emph{Glitch studies manifesto}, Video vortex reader II: Moving
  images beyond YouTube (2011), 336--347.

\bibitem{murdock2017exploration}
Jaimie Murdock, Colin Allen, and Simon DeDeo, \emph{{Exploration and
  exploitation of Victorian science in Darwin's reading notebooks}}, Cognition
  \textbf{159} (2017), 117--126.

\bibitem{thi}
C.~Thi Nguyen, \emph{{The Uses of Aesthetic Testimony}}, The British Journal of
  Aesthetics \textbf{57} (2017), no.~1, 19--36.

\bibitem{nguyen2020games}
\bysame, \emph{Games: Agency as art}, Oxford University Press, USA, 2020.

\bibitem{catarina}
Catarina~Dutilh Novaes, \emph{The dialogical roots of deduction: Historical,
  cognitive, and philosophical perspectives on reasoning}, Cambridge University
  Press, 2020.

\bibitem{putnam}
Hilary Putnam, \emph{Rethinking mathematical necessity}, The New Wittgenstein
  (Alice Crary and Rupert Read, eds.), Routledge, 2000, pp.~245--263.

\bibitem{robson2012aesthetic}
Jon Robson, \emph{Aesthetic testimony}, Philosophy Compass \textbf{7} (2012),
  no.~1, 1--10.

\bibitem{russell}
Bertrand Russell, \emph{Mysticism and logic}, Taylor Garnett Evans \& Co.,
  Watford, Hertsfordshire, UK, 1917.

\bibitem{scholze2022liquid}
Peter Scholze, \emph{Liquid tensor experiment}, Experimental Mathematics
  \textbf{31} (2022), no.~2, 349--354.

\bibitem{tao}
Terrence Tao, \emph{On ``local'' and ``global'' errors in mathematical papers,
  and how to detect them},
  \url{https://terrytao.wordpress.com/advice-on-writing-papers/on-local-and-global-errors-in-mathematical-papers-and-how-to-detect-them/}.

\bibitem{akshay}
Akshay Venkatesh, \emph{Some thoughts on automation and mathematical research},
   (2022), \url{https://www.math.ias.edu/~akshay/research/IASEssay.pdf}.

\bibitem{viteri}
Scott Viteri and Simon DeDeo, \emph{Epistemic phase transitions in mathematical
  proofs}, Cognition \textbf{225} (2022), 105120.

\bibitem{Williams1972}
B.~A.~O. Williams, \emph{Knowledge and reasons}, pp.~1--11, Springer
  Netherlands, Dordrecht, 1972.

\bibitem{zach}
Zachary Wojtowicz and Simon DeDeo, \emph{From probability to consilience: How
  explanatory values implement bayesian reasoning}, Trends in Cognitive
  Sciences \textbf{24} (2020), no.~12, 981--993.

\bibitem{wolpert2023stochastic}
David~H. Wolpert and David~B. Kinney, \emph{Stochastic mathematical systems},
  arXiv \textbf{2209.00543} (2023), math.LO.

\bibitem{walz2015gameful}
Eric Zimmerman, \emph{Manifesto for a ludic century}, The Gameful World:
  Approaches, Issues, Applications, MIT Press, 2015.

\end{thebibliography}
\bibliographystyle{amsplain}

\end{document}